\documentclass[a4paper,11pt]{amsart}

\usepackage{ psfrag, epsfig,amsmath,amssymb,amsthm,array,multicol,verbatim,moreverb,mathrsfs}
\newtheorem{lemma}{Lemma}[section]

\newtheorem{lem}{Lemma}[section]
\numberwithin{equation}{section} \theoremstyle{plain}
\newtheorem{theorem}{Theorem}[section]

\newtheorem{proposition}{Proposition}[section]

\newtheorem{cor}{Corollary}[section]
\newtheorem{rem}{Remark}[section]
\newenvironment{remark}{\begin{rem}\rm}{\smallskip\end{rem}}
\newtheorem*{ex}{Example}
\newenvironment{example}{\begin{ex}\rm}{\smallskip\end{ex}}

\def\esp{\mathbb{E}}
\def\prob{\mathbb{P}}
\def\R{\mathbb{R}}
\def\real{\mathbb{R}}

\def\Z{\mathbb{Z}}

\def\Fcal{\mathcal{F}}

\def\ind{{\rm 1\hspace{-0.90ex}1}}

\def\1{{\bf 1}}

\def\0{{\bf 0}}

\def\Fcal{\mathcal{F}}

\long\def\symbolfootnote[#1]#2{\begingroup
\def\thefootnote{\fnsymbol{footnote}}\footnote[#1]{#2}\endgroup}

\author[Marc Lelarge]{Marc Lelarge}
\title[Tail asymptotics for networks]{Tail asymptotics for monotone-separable networks}
\date{\today}

\begin{document}
\maketitle
\symbolfootnote[0]{{\it MSC 2000 subject
classifications.}  60F10, 
60K25.} 

\symbolfootnote[0]{{\it Key words.} large deviations, queueing
networks.}

\begin{abstract}
A network belongs to the monotone separable class if its state
variables are homogeneous and monotone functions of the epochs of
the arrival process. This framework contains several classical
queueing network models, including generalized Jackson networks,
max-plus networks, polling systems, multiserver queues, and various
classes of stochastic Petri nets. We use comparison relationships
between networks of this class with i.i.d. driving sequences and the
$GI /GI /1/1$ queue to obtain the tail asymptotics of the stationary
maximal dater under light-tailed assumptions for service times. The
exponential rate of decay is given as a function of a logarithmic
moment generating function. We exemplify an explicit computation of
this rate for the case of queues in tandem under various stochastic
assumptions.
\end{abstract}

\section{Introduction}

Consider the $GI/GI/1$ single server queue: we denote
$X_n=\sigma_n-\tau_n$ where $\{\sigma_n\}$ and $\{\tau_n\}$ are
independent and identically distributed (i.i.d) non-negative random
variables, $\sigma_n$ is the amount of service of customer $n$ and
$\tau_n$ is the inter-arrival time between customer $n$ and $n+1$.
Assume that $\esp[X_1]<0$, then the supremum of the random walk
$S_n=X_1+\dots +X_n$ defined by $M:= \sup_{n\geq 1}S_n$ is finite
almost surely and has the same distribution as the stationary
workload of the single server queue. If we assume moreover that
$\esp[\exp(\epsilon X_1)]<\infty$ for some $\epsilon>0$, then the
following asymptotics is standard:
\begin{eqnarray}
\label{casadd}\lim_{x\rightarrow \infty} \frac{1}{x} \log \prob(M>x)
= -\theta^*, \quad \mbox{where $\theta^*=\sup\left\{\theta>0, \:
\log\esp\left[ e^{\theta X}\right]<0\right\}$.}
\end{eqnarray}
Motivated by queueing applications, this case has been extensively
studied in the literature and much finer estimates are available,
see the works of Iglehart \cite{igl} and Pakes \cite{pak}. The main
goal of this paper is to derive analogous results to (\ref{casadd})
for networks.

In the context of a network, we consider the maximal dater $Z$ which
is the time to empty the network when stopping further arrivals.
Clearly in the single server queue, the maximal dater corresponds to
the workload. In the case of queues in tandem, it corresponds to the
end to end delay. Our Theorem \ref{the:main-net} gives the
logarithmic tail asymptotics for the maximal dater of a monotone
separable network. The main difficulty in our task is the absence of
closed form formula for $Z$. The proof of the theorem will proceed
by deriving upper and lower bounds for monotone separable networks.
This class, which was introduced by Baccelli and Foss in
\cite{bacfos:sat}, contains several classical queueing network
models like generalized Jackson networks, max-plus networks, polling
systems and multiserver queues. In this paper, we choose to put a
particular emphasis on tandem queues that fall in the class of open
Jackson networks, and in the class of open (max,plus) systems which
both belong to the class of monotone separable networks. It serves
as a pedagogical example to apply our main theorem under various
stochastic assumptions and it enables us to link our results with
existing asymptotics results from queueing literature.

The paper is structured as follows. In Section \ref{sec:main}, we
give the precise definition of a monotone separable network and its
associated maximal dater. We then give the main result of this paper
in Section \ref{sec:ass}. The case of queues in tandem is dealt with
great details in Section \ref{sec:tand}. In particular, we show that
a kind of phase transition is possible when service times at both
station are dependent. We also link our result to the literature.
Finally technical proofs are deferred to Section \ref{sec:proof}.

\section{Tail asymptotics for monotone-separable
networks}\label{sec:main}

In this paper, we consider open stochastic networks with a single
input process $N$ which is a marked point process with points
$\{T_n\}$ corresponding to exogenous arrival times and marks
$\{\zeta_n\}$ which describe the service times and routing
decisions.

More precisely a stochastic network is described by the following
framework (introduced by Baccelli and Foss \cite{bacfos:sat})
\begin{itemize}
\item The network has a single input point process $N$, with points $\{T_n\}$; for all $m\leq n\in \Z$, let $N_{[m,n]}$ be the $[m,n]$-restriction of $N$, namely the point process with points $\{T_\ell\}_{m\leq \ell\leq n}$.
\item The network has a.s. finite activity for all finite restrictions of $N$: for all $m\leq n\in\Z$, let $X_{[m,n]}(N)$ be the time of last activity in the network, when this one starts empty and is fed by $N_{[m,n]}$. We assume that for all finite $m$ and $n$ as above, $X_{[m,n]}(N)$ is finite.
\end{itemize}
We assume that there exists a set of functions $\{f_\ell\}$,
$f_\ell:\R^\ell\times K^\ell\rightarrow \R$, such that:
\begin{eqnarray}
\label{defX}X_{[m,n]}(N)=f_{n-m+1}\{(T_\ell,\zeta_\ell),\:m\leq
\ell\leq n\},
\end{eqnarray}
for all $n,m$ and $N=\{T_n\}$, where the sequence $\{\zeta_n\}$ is
that describing service times and routing decisions.

\begin{example}
Consider a $G/G/1/\infty  \to ./G/1/\infty$ tandem queue. Denote by
$\{\sigma^{(i)}_n\}$ the sequence of service times in station $i =
1, 2$ and $N=\{T_n\}$ the sequence of arrival times at the first
station. With the notation introduced above, we have
$\zeta_n=(\sigma^{(1)}_n, \sigma^{(2)}_n)$ and the time of last
activity is given by,
\begin{eqnarray}\label{eq:Xtand}
X_{[m,n]}(N)=\sup_{m\leq k\leq n} \left\{T_k + \sup_{k\leq i\leq
n}\sum_{j=k}^i \sigma^{(1)}_j+\sum_{j=i}^n\sigma^{(2)}_j\right\}.
\end{eqnarray}
We refer to the Appendix for an explicit derivation of Equation
(\ref{eq:Xtand}). $X_{[m,n]}(N)$ is simply the last departure time
from the network, when only customers $m,m+1,\dots, n$ enter the
network.
\end{example}

We say that a network described as above is monotone-separable if
the functions $f_n$ are such that the following properties hold for
all input point process $N$:
\begin{enumerate}
\item {\bf Causality:} for all $m\leq n$,
\begin{eqnarray*}
X_{[m,n]}(N)\geq T_n;
\end{eqnarray*}
\item {\bf External monotonicity:} for all $m\leq n$,
\begin{eqnarray*}
X_{[m,n]}(N')\geq X_{[m,n]}(N),
\end{eqnarray*}
whenever $N':=\{T'_n\}$ is such that $T'_n\geq T_n$ for all $n$;
\item {\bf Homogeneity:} for all $c\in \R$ and for all $m\leq n$
\begin{eqnarray*}
X_{[m,n]}(N+c)=X_{[m,n]}(N)+c,
\end{eqnarray*}
where $N+c$ is the point process with points $\{T_n+c\}$;
\item {\bf Separability:} for all $m\leq \ell<n$, if
$X_{[m,\ell]}(N)\leq T_{\ell+1}$, then
\begin{eqnarray*}
X_{[m,n]}(N)=X_{[\ell+1,n]}(N).
\end{eqnarray*}
\end{enumerate}

\begin{remark}
Clearly, tandem queues belong to the class of monotone-separable
networks.
\end{remark}

\subsection{Stability and stationary maximal daters}

In this section, we introduce stochastic assumptions ensuring the
stability of the network. More general results can be found in
Baccelli and Foss \cite{bacfos:sat} and we refer to it for the
statements given in this section without proof.

By definition, for $m\leq n$, the $[m,n]$ maximal dater is
$$Z_{[m,n]}(N):= X_{[m,n]}(N)-T_{n}.$$ Note that
$Z_{[m,n]}(N)$ is a function of $\{\zeta_l\}_{m\leq \ell\leq n}$ and
$\{\tau_l\}_{m\leq \ell\leq n}$ only, where $\tau_n=T_{n+1}-T_n$. In
particular, $Z_n:=Z_{[n,n]}(N)$ is not a function of $N$ (which
makes the notation consistent).

\begin{lemma}\cite{bacfos:sat} {\bf Internal monotonicity of $X$ and $Z$}\\
Under the above conditions, the variables $X_{[m,n]}$ and
$Z_{[m,n]}$ satisfy the internal monotonicity property: for all $N$,
$m\leq n$,
\begin{eqnarray*}
X_{[m-1,n]}(N) &\geq& X_{[m,n]}(N),\\
Z_{[m-1,n]}(N) &\geq& Z_{[m,n]}(N).
\end{eqnarray*}
\end{lemma}
In particular, the sequence $\{Z_{[-n,0]}(N)\}$ is non-decreasing in
$n$. We define the {\it stationary maximal dater} as
\begin{eqnarray*}
Z:= Z_{(-\infty,0]}(N)=\lim_{n\rightarrow\infty}Z_{[-n,0]}(N)\leq
\infty.
\end{eqnarray*}

\begin{example}
In the case of the tandem queues, the stationary maximal dater is
given by:
\begin{eqnarray}
\label{def:Ztand}Z = \sup_{p\leq q\leq 0} \left\{\sum_{k=p}^q
\sigma^{(1)}_k + \sum_{k=q}^0 \sigma^{(2)}_k - (T_0-T_p)\right\},
\end{eqnarray}
and $Z$ is the stationary end to end delay of the network.
\end{example}

\begin{lemma}\label{lem:sub}\cite{bacfos:sat} {\bf Subadditive property of $Z$}\\
Under the above conditions, $\{Z_{[m,n]}(N)\}$ satisfies the
following subadditive property: for all $m\leq \ell<n$, for all $N$,
\begin{eqnarray*}
Z_{[m,n]}(N)\leq Z_{[m,\ell]}(N)+Z_{[\ell+1,n]}(N).
\end{eqnarray*}
\end{lemma}


We assume that the sequence $\{\tau_n,\zeta_n\}_n$ is a sequence of
i.i.d. random variables. The following integrability assumptions are
also assumed to hold (recall that $Z_n=Z_{[n,n]}(N)$ does not depend
on $N$):
\begin{eqnarray*}
\esp[\tau_n]:=a<\infty,\quad \esp[Z_n]<\infty.
\end{eqnarray*}


Denote by $N^0=\{T^0_n\}$ the degenerate input process with
$T^0_n=0$ for all $n$. This degenerate point process plays a crucial
role for the derivation of the stability condition. The following
lemma follows from Lemma \ref{lem:sub} in which we take as
 input point process $N^0$ (note that the constant
$\gamma$ defined below is denoted $\gamma(0)$ in \cite{bacfos:sat}
to emphasize the fact that the input point process is $N^0$).

\begin{lemma}\cite{bacfos:sat}\label{lem:gamma0}
Under the foregoing stochastic assumption, there exists a
non-negative constant $\gamma$ such that
\begin{eqnarray*}
\lim_{n\rightarrow
\infty}\frac{Z_{[-n,0]}(N^0)}{n}=\lim_{n\rightarrow
\infty}\frac{\esp\left[Z_{[-n,0]}(N^0) \right]}{n}=\gamma \: a.s.
\end{eqnarray*}
\end{lemma}

The main result on the stability region is the following:
\begin{theorem}\cite{bacfos:sat}\label{th:stab}
Under the foregoing stochastic assumptions, either $Z=\infty$ a.s.
or $Z<\infty$ a.s.
\begin{enumerate}
\item[(a)] If $\gamma<a$, then $Z<\infty$ a.s.
\item[(b)] If $Z<\infty$ a.s., then $\gamma\leq a$.
\end{enumerate}
\end{theorem}
A proof is given in Section \ref{sec:setupplow}, where we derive an
upper bound and a lower bound that will be used for the study of
large deviations.

\begin{example}
In the case of tandem queues, the constant $\gamma$ is easy to
compute. We have
\begin{eqnarray*}
\lim_{n\to \infty} \sup_{-n\leq q\leq 0}\frac{ \sum_{k=-n}^q
\sigma^{(1)}_k + \sum_{k=q}^0 \sigma^{(2)}_k }{n} = \max\left(
\esp[\sigma^{(1)}_1], \esp[\sigma^{(2)}_1]\right).
\end{eqnarray*}
Hence Theorem \ref{th:stab} gives the standard stability condition:
$\max\left( \esp[\sigma^{(1)}_1],
\esp[\sigma^{(2)}_1]\right)<\esp[\tau_1]$.
\end{example}

\subsection{Moment generating function and tail
asymptotics}\label{sec:ass}


In the rest of the paper, we will make the following assumptions:
\begin{itemize}
\item Assumption {\bf (AA)} on the arrival process into the network $\{T_n\}$:\\
$\{T_n\}$ is a renewal process independent of the service time and
routing sequences $\{\zeta_n\}$.
\item Assumption {\bf (AZ)}: the sequence $\{\zeta_n\}$ is a sequence of i.i.d. random variables, such that the random variable $Z_0$ is light-tailed, i.e. for $\theta$ in a neighborhood of $0$,
\begin{eqnarray*}
\esp[e^{\theta Z_0}]<+\infty.
\end{eqnarray*}
\item Stability: $\gamma<a= \esp[T_1-T_0]$ see Theorem \ref{th:stab}.
\end{itemize}

The subadditive property of $Z$ directly implies the following
property (which is proved in Lemma \ref{lem:lambda}): for any
monotone separable network that satisfies assumption {\bf (AZ)}, the
following limit
\begin{eqnarray}
\label{eq:lambda}\Lambda_Z(\theta) &=& \lim_{n\rightarrow
\infty}\frac{1}{n}\log \esp\left[ e^{\theta Z_{[1,n]}(N^0)}\right],
\end{eqnarray}
exists in $\R\cup\{+\infty\}$ for all $\theta$. Note that the
subadditive property of $Z$ is valid regardless of the point process
$N$ (see Lemma \ref{lem:sub}). Like in the study of the stability of
the network, it turns out that the right quantity to look at is
$Z_{[m,n]}(N^0)$ where $N^0$ is the degenerate input point process
with all its point equal to 0. We also define:
\begin{eqnarray*}
\Lambda_T(\theta) &=& \log\esp\left[ e^{\theta (T_1-T_0)}\right].
\end{eqnarray*}

\begin{theorem}\label{the:main-net}
Under previous assumptions, the tail asymptotics of the stationary
maximal dater is given by,
\begin{eqnarray*}
\lim_{x\rightarrow \infty} \frac{1}{x} \log \prob(Z>x) =
-\theta^*<0,
\end{eqnarray*}
where $\theta^*=\sup \left\{\theta>0,\:
\Lambda_T(-\theta)+\Lambda_Z(\theta)<0 \right\}$.
\end{theorem}

It is relatively easy to see that under our light-tailed assumption
the stationary maximal dater $Z$ will be light-tailed (see Corollary
3 in \cite{bacfos:momtail}). The main contribution of Theorem
\ref{the:main-net} is to give an explicit way of computing the rate
of decay of the tail distribution of $Z$. We refer the interested
reader to \cite{lel:valuetools} for more details on the computation
of $\Lambda_Z$ in the case of (max,plus)-linear networks. In Section
\ref{sec:tand}, we continue the study of our example and deal with
the case of queues in tandem under various stochastic assumptions.
This case of study allows us to show a phase transition phenomena
and to compare our theorem with results of the literature.

Note that in the context of heavy-tailed asymptotics, the moment
generating function is infinite for all $\theta>0$. There is no
general result for the tail asymptotics of the maximal dater of a
monotone separable network. However the methodology derived by
Baccelli and Foss \cite{bacfos:momtail} for subexponential
distributions allows to get exact asymptotics for (max,plus)-linear
networks \cite{bfl} and generalized Jackson networks
\cite{bfl:jack}.

\section{A Case of Study: Queues in tandem}\label{sec:tand}

\subsection{The impact of dependence}

We continue our example and consider a stable $G/G/1/\infty  \to
./G/1/\infty$ tandem queue where $\{\sigma^{(i)}_n\}$ is the
sequence of service times in station $i = 1, 2$ and $\{\tau_n\}$ is
the sequence of inter-arrival times at the first station. We assume
that the sequences $\{(\sigma^{(1)}_n,\sigma^{(2)}_n\}$ and
$\{\tau_n\}_n$ are sequences of i.i.d. random variables such that
$\gamma=\max\left( \esp[\sigma^{(1)}_1],
\esp[\sigma^{(2)}_1]\right)<\esp[\tau_1]$.

We consider two cases:
\begin{itemize}
\item case 1: the sequences $\{\sigma^{(1)}_n\}$,
$\{\sigma^{(2)}_n\}$, $\{\tau_n\}$ are independent.
\item case 2: the sequences $\{\sigma^{(1)}_n\}$ and
$\{\tau_n\}$ are independent and we have
$\sigma^{(2)}_n=\sigma^{(1)}_n$.
\end{itemize}
We denote $\Lambda_i(\theta)= \log\esp[\exp(\theta \sigma^{(i)}_1)]$
and $\delta = \sup\{\theta\geq 0, \: \esp\left[ e^{\theta
\sigma^{(1)}_1}\right]<\infty\}$. A direct application of Theorem
\ref{the:main-net} gives an extension of the results of Ganesh
\cite{ga}:
\begin{cor}\label{cor:tandem}
The tail asymptotics of the stationary end to end delay for two
queues in tandem is given by
\begin{eqnarray*}
\lim_{x\to \infty} \frac{1}{x}\log \prob(Z>x) = -\theta^*,
\end{eqnarray*}
where
\begin{itemize}
\item in case 1: $\theta^*=\min(\theta^1, \theta^2)$ with $\theta^i = \sup\{\theta>0,\:
\Lambda_i(\theta)+\Lambda_T(-\theta)<0\}$;
\item in case 2: $\theta^*=\min\left(\theta^1,\frac{\delta}{2}\right)$.
\end{itemize}
\end{cor}
In case 1, $\theta^i$ is the rate of exponential decay for the tail
distribution of the stationary workload of a single server queue
with interarrival $\tau_n$ and service time $\sigma^{(i)}_n$ and we
have $\theta^*=\min(\theta^1, \theta^2)$. It is well-known that the
stability of such a network is constraint by the "slowest"
component. Here we see that in a large deviations regime, the "bad"
behavior of the network is due to a "bottleneck" component (which is
not necessarily the same as the "slowest" component in average).
Note that in the particular case where the random variables
$\sigma^{(1)}_n, \sigma^{(2)}_n, \tau_n$ are exponentially
distributed with mean $1/\mu^1, 1/\mu^2, a$, we have $\theta^i=
\mu^i-a^{-1}$, and in this case the "slowest" component in average
is also the "bottleneck" component in the large deviations regime.

In the case where the service times are the same at both stations,
Corollary \ref{cor:tandem} shows that the tail behavior of the
random variable $\sigma^{(1)}_1$ described by $\delta$ matters. To
simplify and to get a parametric model, assume that the arrival
process is Poisson with intensity $\lambda:=a^{-1}$ and the service
times are exponentially distributed with mean $1/\mu$. Then
depending on the intensity of the arrival process $\lambda$, two
situations may occur:
\begin{eqnarray*}
\lambda\leq \mu/2 &\Rightarrow & \theta^*= \mu/2,\\
\lambda>\mu/2 &\Rightarrow & \theta^*= \mu-\lambda.
\end{eqnarray*}
In words, we have
\begin{enumerate}
\item if $\lambda<\mu/2$, then the tail asymptotics of the end-to-end
delay is the same as the total service requirement of a single
customer;
\item if $\lambda>\mu/2$, then the tail asymptotics of the
end-to-end delay is the same as in the independent case.
\end{enumerate}
This shows that the behavior of tandems differs from that of a
single server queue. In particular Anantharam \cite{venkat} shows
that for $GI/GI/1$ queues, the build-up of large delays can happen
in one of two ways:
\begin{itemize}
\item If the service times have
exponential tails, then it involves a large number of customers
(whose inter-arrival and service times differ from their mean
values).
\item If the service
times do not have exponential tails, then large delays are caused by
the arrival of a single customer with large service requirement.
\end{itemize}

We see that the first behavior is still valid for queues in tandem
when the service times are independent at each station or if the
intensity of the arrival process is sufficiently large. In contrast,
when the service times are the same at both station, we see that a
single customer can create large delays in the network even under
the assumption of exponential service times (if the intensity of
arrivals is sufficiently small). Note that this phenomena is rather
simple and results intrinsically from the fact that the network
considered is of dimension greater than 2 (i.e. one cannot get such
a phenomena with a single server queue).

\proof Recall that we have
\begin{eqnarray*}
Z_{[1,n]}(N^0) = \sup_{1\leq k\leq n}\sum_{i=1}^k
\sigma^{(1)}_i+\sum_{i=k}^n \sigma^{(2)}_i.
\end{eqnarray*}
In case 1, we have
\begin{eqnarray*} \log \esp\left[ e^{\theta
Z_{[1,n]}(N^0)}\right] &\leq&\log\left(\sum_{k=1}^n
e^{k\Lambda_1(\theta)+(n-k)\Lambda_2(\theta)}\right)
\\
&\leq & \log n +n \max\left(\Lambda_1(\theta),
\Lambda_2(\theta)\right).
\end{eqnarray*}
Hence we have $\Lambda_Z(\theta) =\max\left(\Lambda_1(\theta),
\Lambda_2(\theta)\right)$ and the corollary follows.

In case 2, we have $Z_{[1,n]}(N^0)=\sum_{i=1}^n \sigma^{(1)}_i
+\max_i \sigma^{(1)}_i = \max_i \left(2\sigma^{(1)}_i +\sum_{j\neq
i} \sigma^{(1)}_j\right)$, hence we have
\begin{eqnarray*}
\log\esp\left[e^{\theta Z_{[1,n]}(N^0)} \right] &\geq& \max\left(n
\Lambda_1(\theta), \Lambda_1(2\theta)\right) \mbox{ and,}\\
\log\esp\left[e^{\theta Z_{[1,n]}(N^0)} \right] &\leq&
(n-1)\Lambda_1(\theta) +\log n +\Lambda_1(2\theta).
\end{eqnarray*}
It follows that
\begin{eqnarray*}
\Lambda_Z(\theta) =\left\{\begin{array}{ll}\Lambda_1(\theta)&,\:
\theta <\eta/2\\
\infty&,\: \theta> \eta/2
\end{array} \right.
\end{eqnarray*}
and the corollary follows.
\endproof

\subsection{Comparison with the literature}

In the context of two queues in tandem, if we define
\begin{eqnarray*}
Y_n = \sup_{-n\leq q\leq 0} \sum_{k=-n}^q\sigma^{(1)}_k+\sum_{k=q}^0
\sigma^{(2)}_k - (T_0 - T_{-n}),
\end{eqnarray*}
then we have in view of (\ref{def:Ztand}), $Z=\sup_n Y_n$. The
supremum of a stochastic process has been extensively studied in
queueing theory but we do not know of any general results that would
allow to derive Corollary \ref{cor:tandem}. To end this section and
to make the connection with the existing literature, we state the
following result
\begin{cor}\label{cor}
Consider the system of queues in tandem descried above. Under
assumptions of Theorem \ref{the:main-net} and if
\begin{enumerate}
\item the sequence $\{Y_{n}/n\}$ satisfies a large deviation principle (LDP) with a good rate
function I;
\item there exists $\epsilon>0$ such that
$\Lambda_Z(\theta^*+\epsilon)<\infty$,
\end{enumerate}
where $\theta^*$ is defined as in Theorem \ref{the:main-net}. Then
we have
\begin{eqnarray}
\label{ssq}\lim_{x\rightarrow \infty}\frac{1}{x}\log \prob(Z>x) =
-\theta^* = -\inf_{\alpha>0} \frac{I(\alpha)}{\alpha}.
\end{eqnarray}
\end{cor}

This kind of result has been extensively studied in the queueing
literature (we refer to the work of Duffy, Lewis and Sullivan
\cite{duffy}). However, we see that considering the moment
generating function instead of the rate function allows us to get a
more general result than (\ref{ssq}) since we do not require the
assumption (2) on the tail. Indeed this assumption ensures that the
tail asymptotics of $\prob(Y_{n}>nc)$ for a single $n$ value cannot
dominate those of $\prob(Z>x)$. In this case, equation (\ref{ssq})
has a nice interpretation: the natural drift of the process $Y_{n}$
is $\mu n$, where $\mu<0$. The quantity $I(\alpha)$ can be seen as
the cost for changing the drift of this process to $\alpha>0$. Now
in order to reach level $x$, this drift has to last for a time
$x/\alpha$. Hence the total cost for reaching level $x$ with drift
$\alpha$ is $xI(\alpha)/\alpha$ and the process naturally choose the
drift with the minimal associated cost. As already discussed, this
heuristic is valid only if an assumption as (2) holds. Note also
that in our framework, we do not assume any LDP to hold for the
sequence $\{Y_n/n\}$. In particular, as shown by Corollary
\ref{cor:tandem}, the computation of the moment generating function
$\Lambda_Z$ is much easier than deriving a LDP for $\{Y_n/n\}$.
Lastly, we should stress that for general monotone separable
networks, the maximal dater $Z$ cannot be expressed as the supremum
of a simple stochastic process in which case, the derivation of the
tail asymptotics of $Z$ requires new techniques.

\proof We have only to show that $\theta^* =\inf_{\alpha>0}
\frac{I(\alpha)}{\alpha}$. Thanks to Varadhan's Integral Lemma (see
Theorem 4.3.1 in \cite{demzet}), we have
\begin{eqnarray*}
\lim_{n\rightarrow\infty}\frac{1}{n} \log\esp\left[ e^{\theta
Y_{n}}\right]=\Lambda(\theta)=\sup_{x}\{\theta x-I(x)\},
\end{eqnarray*}
for $\theta< \theta^*+\epsilon$, where
$\Lambda(\theta)=\Lambda_Z(\theta)+\Lambda_T(-\theta)$. Then, the
corollary follows from the following observations for $\theta >0$,
\begin{eqnarray*}
\theta < \inf_{\alpha>0} \frac{I(\alpha)}{\alpha} &\Leftrightarrow&
\theta \alpha -I(\alpha)< 0,\: \forall
\alpha\\
&\Leftrightarrow& \sup_{\alpha }\{\theta \alpha-I(\alpha)\} =
\Lambda(\theta)< 0.
\end{eqnarray*}

\endproof

\section{Proof of the tail asymptotics}\label{sec:proof}
\subsection{Upper $G/G/1/\infty$ queue and lower bound for the
maximal dater}\label{sec:setupplow}

The material of this subsection is not new and may be found in
various references (that are given in what follows). For the sake of
completeness, we include all the proofs. We derive now upper and
lower bounds for the stationary maximal dater $Z$. These bounds
allow to prove Theorem \ref{th:stab} and will be the main tools for
the study of large deviations.

We first derive a lower bound that can also be found in the textbook
\cite{bacbre} see proof of Theorem 2.11.3.
\begin{proposition}\label{prop:low}
We have the following lower bound
\begin{eqnarray*}
Z \geq \sup_{n\geq 0}\left( Z_{[-n,0]}(N^0)+T_{-n}-T_0\right).
\end{eqnarray*}
\end{proposition}
\proof

For $n$ fixed, let $N^n$ be the point process with point
$T^n_j=T_{-n}-T_0$, for all $j$. Then
\begin{eqnarray*}
Z_{[-n,0]}&=& X_{[-n,0]}(N)-T_0 \geq X_{[-n,0]}(N^n)\\
&=& X_{[-n,0]}(N^0)+T_{-n}-T_0=Z_{[-n,0]}(N^0)+T_{-n}-T_0,
\end{eqnarray*}
where we used external monotonicity in the first inequality and
homogeneity between the first and second line.
\endproof

\proof{of Theorem \ref{th:stab} part (b)}

Suppose that $\gamma>a$, then we have
\begin{eqnarray*}
\liminf_{n\rightarrow \infty}\frac{Z_{[-n,0]}(N)}{n}\geq \gamma-a>0,
\end{eqnarray*}
which concludes the proof of part (b).
\endproof

We assume now that $\gamma<a$. We pick an integer $L\geq 1$ such
that
\begin{eqnarray}
\label{defL}\esp\left[Z_{[-L,-1]}(N^0) \right]<La,
\end{eqnarray}
which is possible in view of Lemma \ref{lem:gamma0}. Without loss of
generality, we assume that $T_0=0$. Part (a) of Theorem
\ref{th:stab} follows from the following proposition (that can be
found in \cite{bacfos:momtail}):
\begin{proposition}\label{prop:ub}
The stationary maximal dater $Z$ is bounded from above by the
stationary response time $\hat{R}$ in the $G/G/1/\infty$ queue with
service times
\begin{eqnarray*}
\hat{s}_n:= Z_{[L(n-1)+1,Ln]}(N^0)
\end{eqnarray*}
and inter-arrival times $\hat{\tau}_n:=T_{Ln}-T_{L(n-1)}$, where $L$
is the integer defined in (\ref{defL}). Since
$\esp[\hat{s}_1]<\esp[\hat{\tau}_1]=La$, this queue is stable. With
the convention $\sum_0^{-1}=0$, we have,
\begin{eqnarray*}
Z\leq \hat{s}_0+\sup_{k\geq 0}
\sum_{i=-k}^{-1}\left(\hat{s}_i-\hat{\tau}_{i+1}\right).
\end{eqnarray*}
\end{proposition}
\proof

To an input process $N$, we associate the following upper bound
process, $N^+=\{T_n^+\}$, where $T_n^+=T_{kL}$ if
$n=(k-1)L+1,\dots,kL$. Note that $T_n^+\geq T_n$ for all $n$. Then
for all $n$, since we assumed $T_0=0$, we have thanks to the
external monotonicity,
\begin{eqnarray}
\label{X=Z}X_{[-n,0]}(N)=Z_{[-n,0]}(N)\leq
X_{[-n,0]}(N^+)=Z_{[-n,0]}(N^+).
\end{eqnarray}
We show that for all $k\geq 1$,
\begin{eqnarray}
\label{majZkL}Z_{[-kL+1,0]}(N^+) \leq \hat{s}_0+\sup_{-k+1\leq i\leq
0} \sum_{j=-i}^{-1}(\hat{s}_{j}-\hat{\tau}_{j+1}).
\end{eqnarray}
This inequality will follow from the two next lemmas
\begin{lemma}
\label{lem:ub1}Assume $T_0=0$. For any $m<n\leq 0$,
\begin{eqnarray*}
Z_{[m,0]}(N) \leq Z_{[n,0]}(N)+(Z_{[m,n-1]}(N)-\tau_{n-1})^+.
\end{eqnarray*}
\end{lemma}
\proof

Assume first that $Z_{[m,n-1]}(N)-\tau_{n-1}\leq 0$, which is
exactly $X_{[m,n-1]}(N)\leq T_{n}$. Then by the separability
property, we have
\begin{eqnarray*}
Z_{[m,0]}(N)=X_{[m,0]}(N)=X_{[n,0]}(N)=Z_{[n,0]}(N).
\end{eqnarray*}
Assume now that $Z_{[m,n-1]}(N)-\tau_{n-1}>0$. Let $N'=\{T_j'\}$ be
the input process defined as follows
\begin{eqnarray*}
\forall j\leq n-1,\quad T'_j&=& T_j,\\
\forall j\geq n,\quad T'_j&=& T_j+Z_{[m,n-1]}(N)-\tau_{n-1}.
\end{eqnarray*}
Then we have $N'\geq N$ and $X_{[m,n-1]}(N')\leq T'_{n}$, hence by
the external monotonicity, the separability and the homogeneity
properties, we have
\begin{eqnarray*}
Z_{[m,0]}(N)&=&X_{[m,0]}(N)\leq X_{[m,0]}(N')\\
&=&X_{[n,0]}(N')=X_{[n,0]}(N)+Z_{[m,n-1]}(N)-\tau_{n-1}\\
&=&Z_{[n,0]}(N)+Z_{[m,n-1]}(N)-\tau_{n-1}.
\end{eqnarray*}
\endproof

From this lemma we derive directly
\begin{lemma}\label{lem:Zn0up}
Assume $T_0=0$. For any $n<0$,
\begin{eqnarray*}
Z_{[n,0]}(N) \leq \sup_{n\leq k\leq 0}\left( \sum_{i=k}^{-1}
(Z_i-\tau_{i+1})\right)+Z_0,
\end{eqnarray*}
with the convention $\sum_0^{-1}=0$
\end{lemma}

Applying Lemma \ref{lem:Zn0up} to $Z_{[-kL+1,0]}(N^+)$ gives
(\ref{majZkL}). We now return to the proof of Proposition
\ref{prop:ub}. We have
\begin{eqnarray*}
Z&=& \lim_{k\rightarrow \infty} Z_{[-kL+1,0]}\\
&=& \sup_{k\geq 0}Z_{[-kL+1,0]}(N)\\
&\leq& \sup_{k\geq 0}Z_{[-kL+1,0]}(N^+) \quad \mbox{thanks to (\ref{X=Z})}\\
&\leq& \sup_{k\geq 0}\left( \hat{s}_0+\sup_{-k+1\leq i\leq 0}
\sum_{j=-i}^{-1}(\hat{s}_{j}-\hat{\tau}_{j+1})\right)=\hat{R},\quad
\mbox{thanks to (\ref{majZkL}).}
\end{eqnarray*}
from Lemma \ref{lem:Zn0up}.
\endproof

\subsection{Moment generating function}

\begin{lem}\label{lem:lambda}
The function $\Lambda_Z(.)$ defined by (\ref{eq:lambda}) is a proper
convex function with $\Lambda_Z(\theta)<\infty$ for all
$\theta<\eta$ and $\Lambda_Z(\theta)=\infty$ for all $\theta>\eta$,
where $\eta=\sup\left\{\theta ,\: \esp[\exp(\theta
Z_0)]<\infty\right\}$.
\end{lem}
\proof Let
\begin{eqnarray*}
\label{def:lambdan}\Lambda_{Z,n}(\theta) &=&\log \esp\left[
e^{\theta \frac{Z_{[1,n]}(N^0)}{n}}\right].
\end{eqnarray*}

Thanks to the subadditive property of $Z$, we have,
\begin{eqnarray*}
Z_{[1,n+m]}(N^0) \leq Z_{[1,n]}(N^0)+Z_{[n+1,n+m]}(N^0),
\end{eqnarray*}
and $Z_{[1,n]}(N^0)$ and $Z_{[n+1,n+m]}(N^0)$ are independent. Hence
for $\theta\geq 0$, we have,
\begin{eqnarray*}
\label{sub+}\Lambda_{Z,n+m}((n+m)\theta) \leq
\Lambda_{Z,n}(n\theta)+\Lambda_{Z,m}(m\theta).
\end{eqnarray*}

Hence we can define for any $\theta\geq 0$,
\begin{eqnarray*}
\Lambda_Z(\theta) &=& \lim_{n\rightarrow \infty}\frac{1}{n}\log
\esp\left[ e^{\theta Z_{[1,n]}(N^0)}\right]=\lim_{n\rightarrow
\infty}\frac{\Lambda_{Z,n}(n \theta)}{n}=\inf_{n\geq
1}\frac{\Lambda_{Z,n}(n \theta)}{n},
\end{eqnarray*}
as an extended real number. The fact that $\Lambda_Z$ is a proper
convex function follows from Lemma 2.3.9 of \cite{demzet}. The last
fact follows from,
\begin{eqnarray*}
\Lambda_Z(\theta)\leq\log \esp\left[ e^{\theta Z_{1}}\right] \mbox{
and, } \log \esp\left[ e^{\theta Z_{1}}\right] \leq
\Lambda_{Z,n}(n\theta)\mbox{ for $\theta\geq 0$ and all $n\geq 1$.}
\end{eqnarray*}
\endproof

We define
\begin{eqnarray*}
\Lambda(\theta) = \Lambda_T(-\theta)+\Lambda_Z(\theta) \mbox{ and }
\Lambda_n(\theta) =\Lambda_T(-\theta)+\Lambda_{Z,n}(\theta).
\end{eqnarray*}
Note that $\Lambda_Z(.)$ and $\Lambda_T(.)$ are proper convex
functions, hence $\Lambda(.)$ is a well defined convex function.
Recall that $\theta^*$ is defined as follows:
\begin{eqnarray*}
\theta^*= \sup\{\theta>0,\:\Lambda(\theta)<0\}.
\end{eqnarray*}
The following lemma is used repeatedly in what follows,

\begin{lem}
\label{lemtheta} Under the foregoing assumptions, we have
$\theta^*>0$ and
\begin{eqnarray*}
\Lambda(\theta)<0 && \mbox{if}\quad \theta \in (0,\theta^*),\\
\Lambda(\theta)>0 && \mbox{if}\quad \theta >\theta^*.
\end{eqnarray*}
\end{lem}
\proof Let
\begin{eqnarray}
\label{thetan}\theta_n = \sup\{\theta>0,\: \Lambda_n(n\theta)<0\}.
\end{eqnarray}
We fix $n$ such that $\esp[Z_{[1,n]}(N^0)]\leq na$, which is
possible in view of the stability condition.

We first show that $\theta_n>0$ and
\begin{eqnarray}
\label{ineg1}\Lambda_n(n\theta)<0 && \mbox{if}\quad \theta \in (0,\theta_n),\\
\label{ineg2}\Lambda_n(n\theta)>0 && \mbox{if}\quad \theta>\theta_n
\end{eqnarray}
The function $\theta\mapsto \Lambda_n(n\theta)$ is convex,
continuous and differentiable on $[0,\eta)$. Hence we have
\begin{eqnarray*}
\Lambda_n(n\delta)=\delta
\left(\esp[Z_{[1,n]}(N^0)]-a\right)+o(\delta),
\end{eqnarray*}
which is less than zero for sufficiently small $\delta>0$. Hence,
the set over which the supremum in the definition of $\theta_n$ is
taken is not empty and $\theta_n>0$. Now (\ref{ineg1}) and
(\ref{ineg2}) follow from the definition of $\theta_n$, the
convexity of $\theta\mapsto \Lambda_n(n\theta)$ and the fact that
$\Lambda_n(0)=0$.

We now show that $\theta_n\rightarrow \theta^*$ as $n\rightarrow
\infty$. We have for $\theta\geq 0$
\begin{eqnarray*}
\lim_{n\rightarrow \infty}\frac{\Lambda_n(n\theta)}{n} = \inf_{n\geq
  1} \frac{\Lambda_n(n\theta)}{n} =\Lambda(\theta).
\end{eqnarray*}
Hence for $\theta\geq 0$, we have $\frac{\Lambda_n(n\theta)}{n}\geq
\Lambda(\theta)$ and
\begin{eqnarray*}
\forall \theta\in(0,\theta_n),\quad
\Lambda(\theta)\leq\frac{\Lambda_n(n\theta)}{n}<0.
\end{eqnarray*}
This implies that $\theta^*\geq \theta_n>0$. If $\theta^*<\infty$,
we can choose $\epsilon>0$ such that $\theta^*-\epsilon>0$ and then
we have $\Lambda_n(n(\theta^*-\epsilon))/n\rightarrow
\Lambda(\theta^*-\epsilon)<0$. Hence for sufficiently large $n$, we
have $\frac{\Lambda_n(n(\theta^*-\epsilon))}{n}<0$, hence
$\theta^*-\epsilon\leq\theta_n$, and we proved that
$\theta_n\rightarrow \theta^*$.
$\Lambda(.)$ is a convex function 
and since $\Lambda(0)=0$, the lemma follows.

If $\theta^*=\infty$, we still have $\theta_n\rightarrow \infty$
(that will be needed in proof of Lemma \ref{lem:upper}) by the same
argument as above with $\theta^*-\epsilon$ replaced by any real
number.

\endproof

\subsection{Lower Bound}\label{sec:lower}

\begin{lem}\label{lem:lower}
Under previous assumptions, we have
\begin{eqnarray*}
\liminf_{x\rightarrow \infty} \frac{1}{x} \log \prob(Z>x) \geq
-\theta^*.
\end{eqnarray*}
\end{lem}
\proof We have (see Proposition \ref{prop:low})
\begin{eqnarray}
\label{lowerb}Z \geq \sup_n \left\{
Z_{[-n,0]}(N^0)+T_{-n}-T_0\right\}.
\end{eqnarray}
We denote $Y_n = Z_{[-n,1]}(N^0)+T_{-n}+T_0$, the lemma will follow
from the following fact: $$\liminf_{x\to \infty} \frac{1}{x}\log
\prob(\sup_n Y_n>x)\geq -\theta^*.$$ Note that we have
\begin{eqnarray*}
\lim_{n\to \infty}\frac{1}{n}\log\esp\left[ e^{\theta
Y_n}\right]=\Lambda(\theta).
\end{eqnarray*}
In particular, we are in the setting of G\"artner-Ellis theorem see
Theorem 2.3.6 in \cite{demzet} which will be the main tool of the
proof.

First note that we only need to consider the case $\theta^*<\infty$.
We consider first the case where there exists $\theta> \theta^*$
such that $\Lambda(\theta)<\infty$. 
First note that the function $\theta\mapsto \Lambda(\theta)$ is
convex, hence the left-hand derivatives $\Lambda'(\theta-)$ and the
right-hand derivatives $\Lambda'(\theta+)$ exist for all $\theta>0$.
Moreover, we have $\Lambda'(\theta-)\leq\Lambda'(\theta+)$ and the
function $\theta\mapsto
\frac{1}{2}(\Lambda'(\theta-)+\Lambda'(\theta+))$ is non-decreasing,
hence $\Lambda'(\theta)=\Lambda'(\theta-)=\Lambda'(\theta+)$ except
for $\theta\in \Delta$, where $\Delta$ is at most countable. Since
$\Lambda(\theta)<\infty$ for $\theta>\theta^*$, we have
$\Lambda(\theta^*)=0$ and $\Lambda'(\theta^*+)>0$. To prove this,
assume that $\Lambda'(\theta^*+)=0$. Take $\theta<\theta^*$, thanks
to Lemma \ref{lemtheta}, we have $\Lambda(\theta)<0$. Choose
$\epsilon>0$ such that $0<\Lambda(\theta^*+\epsilon)<\epsilon
|\Lambda(\theta)|$. We have
\begin{eqnarray*}
\frac{\Lambda(\theta^*+\epsilon)}{\epsilon} <
\frac{-\Lambda(\theta)}{\theta^*-\theta},
\end{eqnarray*}
which contradicts the convexity of $\Lambda(\theta)$. Hence, we can
find $t\leq \theta^*+\epsilon$ such that
\begin{eqnarray*}
0<\Lambda(t),\quad t\notin \Delta.
\end{eqnarray*}
Note that these conditions imply $t>\theta^*$ and $\Lambda'(t)\geq
\Lambda'(\theta^*+)>0$.

Thanks to G\"artner-Ellis theorem (Theorem 2.3.6 in \cite{demzet}),
we have
\begin{eqnarray}
\label{eq:gart}\liminf_{n\rightarrow \infty} \frac{1}{n} \log
\prob(Y_{n}>n \alpha) \geq - \inf_{x \in \Fcal,\: x>\alpha}
\Lambda^*(x),
\end{eqnarray}
where $\Fcal$ is the set of exposed point of $\Lambda^*$ and
$\Lambda^*(x) = \sup_{\theta\geq 0}(\theta x-\Lambda(\theta))$. Note
that from the monotonicity of $\theta x-\Lambda(\theta)$ in $x$ as
$\theta$ is fixed, we deduce that $\Lambda^*$ is non-decreasing.
Moreover take $\alpha =\Lambda'(t)$, then $\Lambda^*(\alpha) =
t\alpha-\Lambda(t)$ and $\alpha\in \Fcal$ by Lemma 2.3.9 of
\cite{demzet}.

Given $x>0$, define $n=\lceil x/\alpha\rceil$. We have
\begin{eqnarray*}
\frac{1}{x}\log \prob(\sup_n Y_n>x) \geq
\frac{1}{n\alpha}\log\prob(Y_{n}\geq n\alpha),
\end{eqnarray*}
taking the limit in $x$ and $n$ (while $\alpha=\Lambda'(t)$ is
fixed) gives thanks to (\ref{eq:gart}),
\begin{eqnarray*}
\liminf_{x\rightarrow \infty} \frac{1}{x} \log \prob(\sup_n
Y_n>x)\geq-\frac{t\alpha-\Lambda(t)}{\alpha} \geq -t\geq
-\theta^*-\epsilon.
\end{eqnarray*}

We consider now the case where for all $\theta>\theta^*$, we have
$\Lambda(\theta)=\infty$, i.e. $\theta^*=\eta$ defined in Lemma
\ref{lem:lambda}. Take $K>0$ and define $\tilde{Z}^K_{[n,m]} =
Z_{[n,m]}(N^0)\prod_{i=n}^m\ind(Z_i\leq K)$ and $\tilde{Z}^K =
\sup_{n\geq 0} (\tilde{Z}^K_{[-n,0]}+T_{-n})$. By (\ref{lowerb}), we have $Z\geq
\tilde{Z}^K$. It is easy to see that the proof of Lemma
\ref{lem:lambda} is still valid (note that the subadditive property
carries over to $\tilde{Z}^K_{[n,m]}$) and the following limit exists
\begin{eqnarray*}
\tilde{\Lambda}^K_Z(\theta) = \lim_{n\to \infty}\frac{1}{n}\log\esp\left[
  e^{\theta \tilde{Z}^K_{[1,n]}}\right] = \inf_n\frac{1}{n}\log\esp\left[
  e^{\theta \tilde{Z}^K_{[1,n]}}\right].
\end{eqnarray*}
Moreover thanks to the subadditive property of $Z$, we have
$\prob(\tilde{Z}^K_{[1,n]}\leq nK)=1$, so that
$\tilde{\Lambda}^K_Z(\theta)\leq \theta K$. Hence by the first part of
the proof, we have
\begin{eqnarray*}
\liminf_{x\to\infty}\frac{1}{x}\log\prob(\tilde{Z}^K>x) \geq -\tilde{\theta}^K,
\end{eqnarray*}
with $\tilde{\theta}^K = \sup\{\theta>0,\:
\tilde{\Lambda}^K_Z(\theta)+\Lambda_T(-\theta)<0\}$. We now prove that
$\tilde{\theta}^K\to\eta$ as $K$ tends to infinity which will
conclude the proof.
Note that for any fixed $\theta\geq 0$, the function
$\tilde{\Lambda}^K_Z(\theta)$ is nondecreasing in $K$ and
$\lim_{K\rightarrow \infty}\tilde{\Lambda}^K_Z(\theta) =
\tilde{\Lambda}_Z(\theta)\leq \Lambda_Z(\theta)$. This directly implies
that $\tilde{\theta}^K\geq \eta$.
Take $\theta>\eta$, so that $\Lambda_Z(\theta)=\infty$. If
$\tilde{\Lambda}_Z(\theta)<\infty$, then for all $K$, we have
$\tilde{\Lambda}^K_Z(\theta)\leq\tilde{\Lambda}_Z(\theta)<\infty$. But, we
have $\tilde{\Lambda}^K_Z(\theta)=\inf_n
\frac{1}{n}\log\esp\left[
  e^{\theta \tilde{Z}^K_{[1,n]}}\right]$, so that there exists $n$ such that
\begin{eqnarray*}
\esp\left[ e^{\theta Z_{[1,n]}(N^0)},\:\max(Z_1,\dots ,Z_n)\leq K\right]\leq e^{n(\tilde{\Lambda}^K_Z(\theta)+1)}\leq e^{n(\tilde{\Lambda}_Z(\theta)+1)},
\end{eqnarray*}
but the left-hand side tends to infinity as $K\to \infty$. Hence we
proved that for all $\theta>\eta$, we have
$\tilde{\Lambda}^K_Z(\theta)\to\infty$ as $K\to \infty$. This implies
that $\tilde{\theta}^K\to \eta$ as $K\to \infty$.
\endproof

\subsection{Upper bound}

\begin{lem}\label{lem:upper}
Under previous assumptions, we have
\begin{eqnarray*}
\limsup_{x\rightarrow\infty}\frac{1}{x}\log\prob(Z>x) \leq
-\theta^*.
\end{eqnarray*}
\end{lem}
\proof For $L$ sufficiently large, we have with the convention
$\sum_0^{-1}=0$ (see Proposition \ref{prop:ub}),
\begin{eqnarray*}
Z\leq \sup_{n\geq
0}\left(\sum_{i=-n}^{-1}\hat{s}_i(L)-\hat{\tau}_{i+1}(L)
\right)+\hat{s}_0(L)=:V(L)+\hat{s}_0(L).
\end{eqnarray*}

We will show that under previous assumptions, we have
\begin{eqnarray}
\label{ldupper}\limsup_{x\rightarrow \infty}\frac{1}{x}\log
\prob(V(L)+\hat{s}_0(L)>x) \leq -\theta_L,
\end{eqnarray}
where $\theta_L$ is defined as in (\ref{thetan}) and the lemma will
follow since $\theta_L\rightarrow \theta^*$ as $L$ tends to infinity
(see Lemma \ref{lemtheta}).

First note that for all $\theta\in(0,\theta_L)$, we have
\begin{eqnarray*}
\max\left\{\esp\left[e^{\theta \hat{s}_0(L)} \right],\esp\left[
e^{\theta V(L)}\right]\right\}<\infty.
\end{eqnarray*}
Hence for $\theta\in(0,\theta_L)$, we have $\esp\left[e^{\theta
(V(L)+\hat{s}_0(L))} \right]=\esp\left[e^{\theta V(L)}
\right]\esp\left[e^{\theta \hat{s}_0(L)} \right]\leq A$ for some
finite constant $A$. Hence by Chernoff's inequality,
\begin{eqnarray*}
\prob\left( V(L)+\hat{s}_0(L)\geq x\right)\leq e^{-\theta
x}\esp\left[e^{\theta (V(L)+\hat{s}_0(L))} \right]\leq Ae^{-\theta
x}.
\end{eqnarray*}
Since the above holds for all $0<\theta<\theta_L$, we get
\begin{eqnarray*}
\limsup_{x\rightarrow\infty}\frac{1}{x}\log\prob\left(
V(L)+\hat{s}_0(L)\geq x\right)\leq -\theta_L.
\end{eqnarray*}\endproof

\section*{Acknowledgment}

The author would like to thank Peter Friz for pointing out a mistake
in an earlier version of this paper.

\section*{Appendix: recursion for queues in tandem}
We consider a $G/G/1/\infty  \to ./G/1/\infty$ tandem queue, where
$\{\sigma^{(i)}_n\}$ denotes the sequence of service times in
station $i = 1, 2$ and $N=\{T_n\}$ is the sequence of arrival times
at the first station. For $m\leq k\leq n$, we denote by
$D^{(i)}_{[m,n]}(k)$ the departure time of customer $k$ from station
$i=1,2$ when the network starts empty and is fed by $N_{[m,n]}$.
With the notations introduced in Section \ref{sec:main}, we have
$X_{[m,n]}(N)=D^{(2)}_{[m,n]}(n)$. We now derive the recursion
equations satisfied by the $D_{[m,n]}$'s,
\begin{eqnarray*}
D^{(1)}_{[m,n]}(m) &=& T_m+\sigma^{(1)}_m,\\
D^{(2)}_{[m,n]}(m) &=&
D^{(1)}_{[m,n]}(m)+\sigma^{(2)}_m=T_m+\sigma^{(1)}_m+\sigma^{(2)}_m,\\
D^{(1)}_{[m,n]}(k) &=& \max\left(D^{(1)}_{[m,n]}(k-1), T_k \right)+\sigma^{(1)}_k,\\
D^{(2)}_{[m,n]}(k) &=& \max\left(D^{(2)}_{[m,n]}(k-1),
D^{(1)}_{[m,n]}(k) \right)+\sigma^{(2)}_k,
\end{eqnarray*}
for $m<k\leq n$. From these equations, one can easily check that:
\begin{eqnarray*}
D^{(1)}_{[m,n]}(k) &=& \sup_{m\leq j\leq k}\left\{ T_j +
\sum_{i=j}^k \sigma^{(1)}_i\right\},\\
D^{(2)}_{[m,n]}(k) &=&\sup_{m\leq j\leq k}\left\{ T_j + \sup_{j\leq
\ell\leq k}\sum_{i=j}^\ell \sigma^{(1)}_i+\sum_{i=\ell}^k
\sigma^{(2)}_i\right\},
\end{eqnarray*}
and Equation (\ref{eq:Xtand}) follows.

\bibliographystyle{abbrv}
\bibliography{ex}

\begin{thebibliography}{10}

\bibitem{venkat}
V.~Anantharam.
\newblock How large delays build up in a {$GI/G/1$} queue.
\newblock {\em Queueing Systems Theory Appl.}, 5(4):345--367, 1989.

\bibitem{bacbre}
F.~Baccelli and P.~Br\'emaud.
\newblock {\em Elements of Queueing Theory}.
\newblock Springer-Verlag, 2003.

\bibitem{bacfos:sat}
F.~Baccelli and S.~Foss.
\newblock On the saturation rule for the stability of queues.
\newblock {\em Journal of Applied Probability}, 32:494--507, 1995.

\bibitem{bacfos:momtail}
F.~Baccelli and S.~Foss.
\newblock Moments and tails in monotone-separable stochastic networks.
\newblock {\em Ann. Appl. Probab.}, 14(2):612--650, 2004.

\bibitem{bfl:jack}
F.~Baccelli, S.~Foss, and M.~Lelarge.
\newblock Tails in generalized {J}ackson networks with subexponential
  service-time distributions.
\newblock {\em J. Appl. Probab.}, 42(2):513--530, 2005.

\bibitem{bfl}
F.~Baccelli, M.~Lelarge, and S.~Foss.
\newblock Asymptotics of subexponential max plus networks: the stochastic event
  graph case.
\newblock {\em Queueing Syst.}, 46(1-2):75--96, 2004.

\bibitem{demzet}
A.~Dembo and O.~Zeitouni.
\newblock {\em Large Deviations Techniques and Applications}.
\newblock Springer-Verlag, 1998.

\bibitem{duffy}
K.~Duffy, J.~T. Lewis, and W.~G. Sullivan.
\newblock Logarithmic asymptotics for the supremum of a stochastic process.
\newblock {\em Ann. Appl. Probab.}, 13(2):430--445, 2003.

\bibitem{ga}
A.~Ganesh.
\newblock Large deviations of the sojourn time for queues in series.
\newblock {\em Annals of Operations Research}, 79:3--26, 1998.

\bibitem{igl}
D.~L. Iglehart.
\newblock Extreme values in the {$GI/G/1$} queue.
\newblock {\em Ann. Math. Statist.}, 43:627--635, 1972.

\bibitem{lel:valuetools}
M.~Lelarge.
\newblock Tail asymptotics for discrete event systems.
\newblock {\em VALUETOOLS}, 2006.

\bibitem{pak}
A.~G. Pakes.
\newblock On the tails of waiting-time distributions.
\newblock {\em J. Appl. Probability}, 12(3):555--564, 1975.

\end{thebibliography}
Marc Lelarge\footnote{This work was partially supported by Science
Foundation Ireland Research Grant No. SFI 04/RP1/I512.}
\\ ENS-INRIA\\
45 rue d'Ulm\\
75005 Paris, France
\\  {\tt e-mail : marc.lelarge@ens.fr}
\end{document}